\documentclass[11pt]{article}

\usepackage{latexsym,color,graphicx,amsmath,amssymb,amsthm}
\usepackage{url}

\newtheorem{theorem}{Theorem}
\newtheorem{lemma}[theorem]{Lemma}

\newtheorem{corollary}[theorem]{Corollary}

\def\reals{{\mathbb R}}

\def\eps{{\varepsilon}}

\def\cplx{{\mathbb C}}

\def\P{{\mathbb P}}

\def\deg{{\mathsf{deg}}}

\newcommand{\Deg}{D} 
\makeatletter
\newcommand{\ProofEndBox}{{\ifhmode\unskip\nobreak\hfil\penalty50 \else
          \leavevmode\fi\quad\vadjust{}\nobreak\hfill$\Box$
            \finalhyphendemerits=0 \par}}
\makeatother

\newcommand{\proofend}{\ProofEndBox\smallskip}

\begin{document}

\title{Incidences between points and lines in three dimensions\thanks{%
Work on this paper by Noam Solomon and Micha Sharir was
supported by Grant 892/13 from the Israel Science Foundation.
Work by Micha Sharir was also supported
by Grant 2012/229 from the U.S.--Israel Binational Science Foundation,
by the Israeli Centers of Research Excellence (I-CORE)
program (Center No.~4/11), and
by the Hermann Minkowski-MINERVA Center for Geometry
at Tel Aviv University.
}}

\author{
Micha Sharir\thanks{%
School of Computer Science, Tel Aviv University,
Tel Aviv 69978, Israel.
{\sl michas@post.tau.ac.il} }
\and
Noam Solomon\thanks{%
School of Computer Science, Tel Aviv University,
Tel Aviv 69978, Israel.
{\sl noam.solom@gmail.com} }
}

\maketitle

\begin{abstract}
We give a fairly elementary and simple proof that shows that the number of incidences
between $m$ points and $n$ lines in $\reals^3$, so that no plane contains more than $s$
lines, is
$$
O\left(m^{1/2}n^{3/4}+ m^{2/3}n^{1/3}s^{1/3} + m + n\right)
$$
(in the precise statement, the constant of proportionality of the first and third terms
depends, in a rather weak manner, on the relation between $m$ and $n$).

This bound, originally obtained by Guth and Katz~\cite{GK2} as a major step in their solution
of Erd{\H o}s's distinct distances problem, is also a major new result in incidence geometry,
an area that has picked up considerable momentum in the past six years. Its original proof
uses fairly involved machinery from algebraic and differential geometry, so it is highly
desirable to simplify the proof, in the interest of better understanding the geometric
structure of the problem, and providing new tools for tackling similar problems. This has
recently been undertaken by Guth~\cite{Gu14}. The present paper presents a different and
simpler derivation, with better bounds than those in \cite{Gu14}, and without the restrictive
assumptions made there. Our result has a potential for
applications to other incidence problems in higher dimensions.
\end{abstract}


\section{Introduction}

Let $P$ be a set of $m$ distinct points in $\reals^3$ and let $L$ be a set of
$n$ distinct lines in $\reals^3$. Let $I(P,L)$ denote the number of
incidences between the points of $P$ and the lines of $L$; that is,
the number of pairs $(p,\ell)$ with $p\in P$, $\ell\in L$, and
$p\in\ell$. If all the points of $P$ and all the lines of $L$ lie in
a common plane, then the classical Szemer\'edi--Trotter
theorem~\cite{SzT} yields the worst-case tight bound
\begin{equation} \label{inc2}
I(P,L) = O\left(m^{2/3}n^{2/3} + m + n \right) .
\end{equation}
This bound clearly also holds in three dimensions, by
projecting the given lines and points onto some generic plane.
Moreover, the bound will continue to be worst-case tight by placing
all the points and lines in a common plane, in a configuration that
yields the planar lower bound.

In the 2010 groundbreaking paper of Guth and Katz~\cite{GK2}, an improved
bound has been derived for $I(P,L)$, for a set $P$ of $m$ points and a set
$L$ of $n$ lines in $\reals^3$, provided that not too many lines of $L$ lie
in a common plane. Specifically, they showed:\footnote{%
  We skip over certain subtleties in their bound: They also assume that no
  \emph{regulus} contains more than $s$ input lines, but then they are able
  also to bound the number of intersection points of the lines. Moreover,
  if one also assumes that each point is incident to at least three lines
  then the term $m$ in the bound can be dropped.}

\begin{theorem}[Guth and Katz~\cite{GK2}]
\label {ttt}
Let $P$ be a set of $m$ distinct points and $L$ a set of $n$ distinct lines
in $\reals^3$, and let $s\le n$ be a parameter,
such that
no plane contains more than $s$ lines of $L$. Then
$$
I(P,L) = O\left(m^{1/2}n^{3/4} + m^{2/3}n^{1/3}s^{1/3} + m + n\right).
$$
\end{theorem}
This bound was a major step in the derivation of the main result of \cite{GK2}, which was
to prove an almost-linear lower bound on the number of distinct distances determined by any
finite set of points in the plane, a classical problem posed by Erd{\H o}s in 1946~\cite{Er46}.
Their proof uses several nontrivial tools from algebraic and differential geometry, most
notably the Cayley--Salmon theorem on osculating lines to algebraic surfaces in $\reals^3$,
and additional properties of ruled surfaces. All this machinery comes on top of the main
innovation of Guth and Katz, the introduction of the \emph{polynomial partitioning technique};
see below.

In this paper, we provide a simple derivation of this bound, which bypasses most of the
techniques from algebraic geometry that are used in the original proof. A recent related study
by Guth~\cite{Gu14} provides another simpler derivation of a similar bound, but (a) the bound
obtained in \cite{Gu14} is slightly worse, involving extra factors of the form $m^\eps$, for
any $\eps>0$, and (b) the assumptions there are stronger, namely that no algebraic surface
of degree at most $c_\eps$, a (potentially large) constant that depends on $\eps$, contains
more than $s$ lines of $L$ (in fact, Guth considers in \cite{Gu14} only the case $s=\sqrt{n}$).
It should be noted, though, that Guth also manages to derive a (slightly weaker but still)
near-linear lower bound on the number of distinct distances.

As in the classical work of Guth and Katz~\cite{GK2}, and in the follow-up study of Guth~\cite{Gu14},
here too we use the polynomial partitioning method, as pioneered in \cite{GK2}.
The main difference between our approach and those of \cite{Gu14,GK2} is the choice of the degree
of the partitioning polynomial. Whereas Guth and Katz \cite{GK2} choose a large degree,
and Guth~\cite{Gu14} chooses a constant degree, we choose an intermediate degree.
This reaps many benefits from both the high-degree and the constant-degree approaches,
and pays a small price in the bound (albeit much better than in \cite{Gu14}). Specifically,
our main result is a simple and fairly elementary derivation of the following result.

\begin{theorem}
\label {th:main}
Let $P$ be a set of $m$ distinct points and $L$ a set of $n$ distinct lines
in $\reals^3$, and let $s\le n$ be a parameter,
such that
no plane contains more than $s$ lines of $L$. Then
\begin{equation} \label{st}
I(P,L) \le A_{m,n} \left(m^{1/2}n^{3/4} + m\right) + B\left( m^{2/3}n^{1/3}s^{1/3} + n\right) ,
\end{equation}
where $B$ is an absolute constant, and, for another suitable absolute constant $b>1$,
\begin{equation} \label{amn}
A_{m,n} =
O\left( b^{\frac{\log (m^2n)}{\log (n^3/m^2)}} \right) ,\quad
\text{for $m \le n^{3/2}$} , \quad\text{and}\quad\
O\left( b^{\frac{\log (m^3/n^4)}{\log (m^2/n^3)}} \right) ,\quad
\text{for $m \ge n^{3/2}$} .
\end{equation}
\end{theorem}

\noindent{\bf Remarks.}
(1) Only the range $\sqrt{n} \le m \le n^2$ is of interest; outside this range,
regardless of the dimension of the ambient space, we have the well known and trivial
upper bound $O(m+n)$.

\smallskip

\noindent
(2) The term $m^{2/3}n^{1/3}s^{1/3}$ comes from the planar Szemer\'edi--Trotter bound
(\ref{inc2}), and is unavoidable, as it can be attained if we densely ``pack'' points
and lines into planes, in patterns that realize the bound in (\ref{inc2}).

\smallskip

\noindent
(3) Ignoring this term, the two terms $m^{1/2}n^{3/4}$ and $m$ ``compete'' for dominance;
the former dominates when $m\le n^{3/2}$ and the latter when $m\ge n^{3/2}$. Thus the
bound in (\ref{st}) is qualitatively different within these two ranges.

\smallskip

\noindent
(4) The threshold $m=n^{3/2}$ also arises in the related
problem of \emph{joints} (points incident to at least three non-coplanar lines) in a
set of $n$ lines in 3-space; see \cite{GK}.

A concise rephrasing of the bound in (\ref{st}) and (\ref{amn}) is as follows.
We partition each of the ranges $m\le n^{3/2}$, $m > n^{3/2}$
into a sequence of subranges $n^{\alpha_{j-1}} < m \le n^{\alpha_j}$, $j=0,1,\ldots$
(for $m\le n^{3/2}$), or $n^{\alpha_{j-1}} > m \ge n^{\alpha_j}$, $j=0,1,\ldots$
(for $m \ge n^{3/2}$), so that within each range the bound asserted in the theorem holds
for some fixed constant of proportionality (denoted as $A_{m,n}$ in the bound),
where these constants vary with $j$, and grow, exponentially in $j$,
as prescribed in (\ref{amn}), as $m$ approaches $n^{3/2}$
(from either side). Informally, if we keep $m$ ``sufficiently away'' from $n^{3/2}$,
the bound in (\ref{st}) holds with a fixed constant of proportionality.
Handling the ``border range'' $m\approx n^{3/2}$ is also fairly straightforward, although,
to bypass the exponential growth of the constant of proportionality, it results in a slightly
different bound; see below for details.

Our proof is elementary to the extent that, among other things, it avoids any
explicit handling of \emph{singular} and \emph{flat} points on the zero set of the
partitioning polynomial. While these notions are relatively easy to handle in three
dimensions (see, e.g., \cite{EKS,GK}), they become more complex notions in higher
dimensions (as witnessed, for example, in our companion work on the four-dimensional
setting~\cite{SS4d}), making proofs based on them harder to extend.

Additional merits and features of our analysis are discussed in detail in the 
concluding section.  In a nutshell, the main merits are: 

\noindent{\bf (i)} 
We use two separate partitioning polynomials. The first one is of ``high'' degree,
and is used to prune away some points and lines, and to establish useful properties 
of the surviving points and lines. The second partitioning step, using a polynomial
of ``low'' degree, is then applied, from scratch, to the surviving input, exploiting
the properties established in the first step. This idea seems to have a potential
for further applications.

\noindent{\bf (ii)} 
Because of the way we use the polynomial partitioning technique, we need induction
to handle incidences within the cells of the second partition. One of the nontrivial
achievements of our technique is the ability to retain The ``planar'' term 
$O(m^{2/3}n^{1/3}s^{1/3})$ in the bound in (\ref{st}) through the inductive process.
Without such care, this term does not ``pass well'' through the induction, which has
been a sore issue in several recent works on related problems (see \cite{SSS,SSZ,surf-socg}).
This is one of the main reasons for using two separate partitioning steps.


\paragraph{Background.}
Incidence problems have been a major topic in combinatorial and
computational geometry for the past thirty years, starting with
the aforementioned Szemer\'edi-Trotter bound \cite{SzT} back in 1983.
Several techniques, interesting in their own right, have been
developed, or adapted, for the analysis of incidences, including
the crossing-lemma technique of Sz\'ekely~\cite{Sz}, and the use of
cuttings as a divide-and-conquer mechanism (e.g., see~\cite{CEGSW}).
Connections with range searching and related algorithmic problems in computational
geometry have also been noted, and studies of the Kakeya problem
(see, e.g., \cite{T}) indicate the connection between this problem and
incidence problems. See Pach and Sharir~\cite{PS} for a
comprehensive (albeit a bit outdated) survey of the topic.

The landscape of incidence geometry has dramatically changed in the
past six years, due to the infusion, in two groundbreaking papers
by Guth and Katz~\cite{GK,GK2}, of new tools and techniques drawn
from algebraic geometry. Although their two direct goals have been to obtain
a tight upper bound on the number of joints in a set of lines in three
dimensions \cite{GK}, and a near-linear lower bound for the classical distinct
distances problem of Erd{\H o}s \cite{GK2}, the new tools have quickly
been recognized as useful for incidence bounds.
See \cite{EKS,KMSS,KMS,SSZ,SoTa,Za1,Za2} for a sample of recent
works on incidence problems that use the new algebraic machinery.

The simplest instances of incidence problems involve points and
lines, tackled by Szemer\'edi and Trotter in the plane~\cite{SzT},
and by Guth and Katz in three dimensions~\cite{GK2}. Other recent
studies on incidence problems include incidences between points and
lines in four dimensions (Sharir and Solomon~\cite{surf-socg,SS4d}), and
incidences between points and circles in three dimensions (Sharir,
Sheffer and Zahl~\cite{SSZ}), not to mention incidences with
higher-dimensional surfaces, such as in \cite{BS,KMSS,SoTa,Za1,Za2}.
In a companion paper (with Sheffer)~\cite{SSS}, we study the general case of
incidences between points and curves in any dimension, and
derive reasonably sharp bounds (albeit weaker in several respects
than the one derived here).

That tools from algebraic geometry form the major key for successful
solution of difficult problems in combinatorial geometry, came as
a big surprise to the community. It has lead to intensive research of the new tools,
aiming to extend them and to find new applications. A major purpose of this study,
as well as of Guth~\cite{Gu14}, is to show that one can still tackle successfully
the problems using less heavy algebraic machinery. This offers a new, simplified,
and more elementary approach, which we expect to prove potent for other applications too,
such as those just mentioned.
Looking for simpler, yet effective techniques that would be easier to extend to
more involved contexts (such as incidences in higher dimensions) has been our
main motivation for this study.

A more detailed supplementary discussion (which would be premature at this point) of the
merits and other issues related to our technique is given in a concluding section.



\section{Proof of Theorem~\ref{th:main}} \label{sec:pf1}

The proof proceeds by induction on $m$. As already mentioned, the bound in (\ref{st})
is qualitatively different in the two ranges $m\le n^{3/2}$ and $m\ge n^{3/2}$. The
analysis bifurcates accordingly. While the general flow is fairly similar in both
cases, there are many differences too.

\paragraph{The case $m < n^{3/2}$.}
We partition this range into a sequence of ranges $m\le n^{\alpha_0}$,
$n^{\alpha_0} < m \le n^{\alpha_1},\ldots$, where $\alpha_0 = 1/2$
and the sequence $\{\alpha_j\}_{j\ge 0}$ is increasing and converges
to $3/2$. More precisely, as our analysis will show, we can take
$\alpha_j = \frac32 - \frac{2}{j+2}$, for $j\ge 0$.
The induction is actually on the index $j$ of the range
$n^{\alpha_{j-1}}<m \le n^{\alpha_j}$, and establishes (\ref{st})
for $m$ in this range, with a coefficient $A_j$ (written in
(\ref{st}, \ref{amn}) as $A_{m,n}$) that increases with $j$.
This paradigm has already been used in Sharir et al.~\cite{SSZ}
and in Zahl~\cite{Za2}, for related incidence problems, albeit in
a somewhat less effective manner; see the discussion at the end of the paper.

The base range of the induction is $m\le \sqrt n$, where the
trivial general upper bound on point-line incidences, in any
dimension, yields $I=O(m^2+n)=O(n)$, so (\ref{st}) holds for a
sufficiently large choice of the initial constant $A_0$.

Assume then that (\ref{st}) holds for all $m \le n^{\alpha_{j-1}}$
for some $j \ge 1$, and consider an instance of the problem with
$n^{\alpha_{j-1}} < m \le n^{3/2}$ (the analysis will force us
to constrain this upper bound in order to complete the induction step,
thereby obtaining the next exponent $\alpha_j$).

Fix a parameter $r$, whose precise value will be chosen later
(in fact, and this is a major novelty of our approach, there will
be two different choices for $r$---see below), and apply the
polynomial partitioning theorem of Guth and Katz (see \cite{GK2} and
\cite[Theorem 2.6]{KMS}), to obtain an $r$-partitioning trivariate (real)
polynomial $f$ of degree $D=O(r^{1/3})$. That is, every connected
component of $\reals^3\setminus Z(f)$ contains at most $m/r$ points of $P$,
where $Z(f)$ denotes the zero set of $f$. By Warren's theorem~\cite{w-lbanm-68}
(see also \cite{KMS}), the number of components of $\reals^3\setminus Z(f)$ is
$O(D^3) = O(r)$.

Set $P_1:= P\cap Z(f)$ and $P'_1:=P\setminus P_1$. A major recurring theme
in this approach is that,
although the points of $P'_1$ are more or less evenly partitioned
among the cells of the partition, no nontrivial bound can be
provided for the size of $P_1$; in the worst case, all the points of
$P$ could lie in $Z(f)$.  Each line $\ell\in L$ is either fully
contained in $Z(f)$ or intersects it in at most $D$ points (since
the restriction of $f$ to $\ell$ is a univariate polynomial of degree
at most $D$). Let $L_1$ denote the subset of lines of $L$ that are
fully contained in $Z(f)$ and put $L'_1 = L\setminus L_1$. We then have
$$
I(P,L) = I(P_1,L_1) + I(P_1,L'_1) + I(P'_1,L'_1) .
$$
We first bound $I(P_1,L'_1)$ and $I(P'_1,L'_1)$. As already observed, we have
$$
I(P_1,L'_1) \le |L'_1|\cdot D \le nD .
$$
We estimate $I(P'_1,L'_1)$ as follows. For each (open) cell $\tau$ of
$\reals^3\setminus Z(f)$, put $P_\tau = P\cap \tau$ (that is,
$P'_1\cap\tau$), and let $L_\tau$ denote the set of the lines
of $L'_1$ that cross $\tau$; put $m_\tau = |P_\tau| \le m/r$, and
$n_\tau = |L_\tau|$. Since every line $\ell\in L'_1$ crosses at
most $1+D$ components of $\reals^3\setminus Z(f)$, we have
$$
\sum_\tau n_\tau \le n(1+D) ,\quad\quad\text{and}\quad\quad
I(P'_1,L'_1) = \sum_{\tau} I(P_\tau,L_\tau).
$$
For each $\tau$ we use the trivial bound
$I(P_\tau,L_\tau) = O(m_\tau^2+n_\tau)$. Summing over the cells, we get
$$
I(P'_1,L'_1) = \sum_{\tau} I(P_\tau,L_\tau) =
O\left(r\cdot(m/r)^2 + \sum_\tau n_\tau \right)
= O\left(m^2/r + nD \right) = O(m^2/D^3 + nD) .
$$
For the initial value of $D$, we take $D = m^{1/2}/n^{1/4}$ (which we
get from a suitable value of $r=\Theta(D^3)$), and get the bound
$$
I(P'_1,L'_1) + I(P_1,L'_1) = O(m^{1/2}n^{3/4}) .
$$
This choice of $D$ is the one made in \cite{GK2}. It is sufficiently large
to control the situation in the cells, by the bound just obtained, but requires
heavy-duty machinery from algebraic geometry to handle the situation on $Z(f)$.

We now turn to $Z(f)$, where we need to estimate $I(P_1,L_1)$. Since all the incidences involving
any point in $P'_1$ and/or any line in $L'_1$ have already been accounted for, we
discard these sets, and remain with $P_1$ and $L_1$ only. We ``forget''
the preceding polynomial partitioning step, and start afresh, applying a
new polynomial partitioning to $P_1$ with a polynomial $g$ of degree $E$,
which will typically be much smaller than $D$, but still non-constant.

Before doing this, we note that the set of lines $L_1$ has a special
structure, because all its lines lie on the algebraic surface $Z(f)$,
which has degree $D$. We exploit this to derive the following lemmas.
We emphasize, since this will be important later on in the analysis, that
Lemmas~\ref{donpl}--\ref{incone} hold for any choice of ($r$ and) $D$.

We note that in general the partitioning polynomial $f$ may be reducible,
and apply some of the following arguments to each irreducible factor separately.
Clearly, there are at most $D$ such factors.

\begin{lemma} \label{donpl}
Let $\pi$ be a plane which is not a component of $Z(f)$.
Then $\pi$ contains at most $D$ lines of $L_1$.
\end{lemma}
\noindent{\bf Proof.}
Suppose to the contrary that $\pi$ contains at least $D+1$ lines of $L$.
Every generic line $\lambda$ in $\pi$ intersects these lines in at least
$D+1$ distinct points, all belonging to $Z(f)$. Hence $f$ must vanish
identically on $\lambda$, and it follows that $f\equiv 0$ on $\pi$, so
$\pi$ is a component of $Z(f)$, contrary to assumption.
\proofend

\begin{lemma} \label{incpl}
The number of incidences between the points of $P_1$ that lie in the
planar components of $Z(f)$ and the lines of $L_1$, is
$O(m^{2/3}n^{1/3}s^{1/3} + nD)$.
\end{lemma}
\noindent{\bf Proof.}
Clearly, $f$ can have at most $D$ linear factors, and thus $Z(f)$ can contain
at most $D$ planar components. Enumerate them as $\pi_1,\ldots,\pi_k$, where
$k\le D$. Let $\tilde{P}_1$ denote the subset of the points of $P_1$ that lie
in these planar components. Assign each point of $\tilde{P}_1$ to
the first plane $\pi_i$, in this order, that contains it, and assign each line of $L_1$
to the first plane that fully contains it; some lines might not be assigned at
all in this manner. For $i=1,\ldots,k$, let $\tilde{P}_i$ denote the set of
points assigned to $\pi_i$, and let $\tilde{L}_i$ denote the set of lines
assigned to $\pi_i$. Put $m_i = |\tilde{P}_i|$ and $n_i = |\tilde{L}_i|$.
Then $\sum_i m_i \le m$ and $\sum_i n_i \le n$; by assumption, we also have
$n_i\le s$ for each $i$. Then
$$
I(\tilde{P}_i,\tilde{L}_i) = O(m_i^{2/3}n_i^{2/3} + m_i + n_i)
= O(m_i^{2/3}n_i^{1/3}s^{1/3} + m_i + n_i) .
$$
Summing over the $k$ planes, we get, using H\"older's inequality,
\begin{align*}
\sum_i & I(\tilde{P}_i,\tilde{L}_i) = \sum_i O(m_i^{2/3}n_i^{1/3}s^{1/3} + m_i + n_i) \\
& = O\left( \left( \sum_i m_i \right)^{2/3} \left( \sum_i n_i \right)^{1/3}s^{1/3} + m + n\right)
= O\left( m^{2/3}n^{1/3}s^{1/3} + m + n\right) .
\end{align*}
We also need to include incidences between points
$p\in \tilde{P}_1$ and lines $\ell\in L_1$ not assigned to the same plane as $p$
(or not assigned to any plane at all).
Any such incidence $(p,\ell)$ can be charged (uniquely) to the intersection point
of $\ell$ with the plane $\pi_i$ to which $p$ has been assigned. The number of such
intersections is $O(nD)$, and the lemma follows.
\proofend

\begin{lemma} \label{md2}
Each point $p\in Z(f)$ is incident to at most $D^2$ lines of $L_1$,
unless $Z(f)$ has an irreducible component that is either a plane containing
$p$ or a cone with apex $p$.
\end{lemma}
\noindent{\bf Proof.}
Fix any line $\ell$ that passes through $p$, and write its parametric
equation as $\{p+tv \mid t\in\reals\}$, where $v$ is the direction
of $\ell$. Consider the Taylor expansion of $f$ at $p$ along $\ell$
$$
f(p+tv) = \sum_{i=1}^D \frac{1}{i!}F_i(p;v)t^i ,
$$
where $F_i(p;v)$ is the $i$-th order derivative of $f$ at $p$ in direction $v$;
it is a homogeneous polynomial in $v$ ($p$ is considered fixed) of degree $i$,
for $i=1,\ldots,D$. For each line $\ell\in L_1$ that passes through $p$, $f$
vanishes identically on $\ell$, so we have $F_i(p;v)=0$ for each $i$. Assuming
that $p$ is incident to more than $D^2$ lines of $L_1$, we conclude that the
homogeneous system
\begin{equation} \label{allf0}
F_1(p;v) = F_2(p;v) = \cdots = F_D(p;v) = 0
\end{equation}
has more than $D^2$ (projectively distinct) roots. The classical
B\'ezout's theorem, applied in the projective plane where the
directions $v$ are represented (e.g., see~\cite{CLO}), asserts that,
since all these polynomials are of degree at most $D$, each pair of
polynomials $F_i(p;v)$, $F_j(p;v)$ must have a common factor. The
following slightly more involved inductive argument shows that in
fact all these polynomials must have a common factor.\footnote{%
  See also~\cite{RSZ} for a similar observation.}
\begin {lemma}
Let $f_1,\ldots, f_n \in \cplx[x,y,z]$ be $n$ homogeneous polynomials
of degree at most $D$. If $|Z(f_1,\ldots, f_n)|>D^2$, then all the
$f_i$'s have a nontrivial common factor.
\end {lemma}

\noindent{\bf Proof.}
The proof is via induction on $n$. The case $n=2$ is precisely the
classical B\'ezout's theorem in the projective plane. Assume that the
inductive claim holds for $n-1$ polynomials. By assumption,
$|Z(f_1,\ldots, f_{n-1})|\ge |Z(f_1,\ldots, f_{n})|>D^2$,
so the induction hypothesis implies that there is
a polynomial $g$ that divides $f_i$, for $i=1,\ldots, n-1$; assume,
as we may, that $g=GCD(f_1,\ldots, f_{n-1})$. If there are more than
$\deg(g)\deg(f_{n})$ points in $Z(g,f_{n})$, then again, by the classical
B\'ezout's theorem in the projective plane, $g$ and $f_n$ have a nontrivial
common factor, which is then also a common factor of
$f_i$, for $i=1,\ldots, n$, completing the proof. Otherwise, put
$\tilde f_i = f_i/g$, for $i=1,\ldots, n-1$. Notice that
$Z(f_1,\ldots, f_{n-1}) = Z(\tilde f_1,\ldots, \tilde f_{n-1})\cup Z(g)$,
implying that each point of $Z(f_1,\ldots, f_n)$ belongs either to
$Z(g) \cap Z(f_n)$ or to $Z(\tilde f_1,\ldots,\tilde f_{n-1})\cap Z(f_n)$.
As $|Z(f_1,\ldots, f_n)|>D^2$ and
$|Z(g,f_{n})| \le \deg(g)\deg(f_n)\le \deg(g)D$, it follows that
$$
|Z(\tilde f_1,\ldots, \tilde f_{n-1})| \ge
|Z(\tilde f_1,\ldots, \tilde f_{n-1}, f_n)|
\ge (D-\deg(g))D>(D-\deg(g))^2.
$$
Hence, applying the induction hypothesis to the polynomials
$\tilde f_1,\ldots, \tilde f_{n-1}$ (all of degree at most $D-\deg(g)$),
we conclude that they have a nontrivial common factor, contradicting the fact
that $g$ is the greatest common divisor of $f_1,\ldots, f_{n-1}$. \proofend

Continuing with the proof of Lemma~\ref{md2},
there is an infinity of directions $v$ that satisfy (\ref{allf0}), so there
is an infinity of lines passing through $v$ and contained in $Z(f)$. The union of these
lines can be shown to be a two-dimensional algebraic variety,\footnote{%
  It is simply the variety given by the equations (\ref{allf0}), rewritten as
  $F_1(p;x-p) = F_2(p;x-p) = \cdots = F_D(p;x-p) = 0$. It is two-dimensional
  because it is contained in $Z(f)$, hence at most two-dimensional, and it cannot
  be one-dimensional since it would then consist of only finitely many lines
  (see, e.g., \cite[Lemma 2.3]{SS4d}).}
contained in $Z(f)$, so $Z(f)$ has an irreducible component that is either a
plane through $p$ or a cone with apex $p$, as claimed.
\proofend

\begin{lemma} \label{incone}
The number of incidences between the points of $P_1$ that lie in the (non-planar)
conic components of $Z(f)$, and the lines of $L_1$, is $O(m + nD)$.
\end{lemma}
\noindent{\bf Proof.}
Let $\sigma$ be such an (irreducible) conic component of $Z(f)$ and let $p$
be its apex. We observe that $\sigma$ cannot contain any line that is not incident to
$p$, because such a line would span with $p$ a plane contained in $\sigma$,
contradicting the assumption that $\sigma$ is irreducible and non-planar.
It follows that the number of incidences between $P_\sigma:=P_1\cap\sigma$ and
$L_\sigma$, consisting of the lines of $L_1$ contained in $\sigma$, is thus
$O(|P_\sigma|+|L_\sigma|)$ ($p$ contributes $|L_\sigma|$ incidences, and every
other point at most one incidence). Applying a similar ``first-come-first-serve''
assignment of points and lines to the conic components of $Z(f)$, as we did for
the planar components in the proof of lemma~\ref{incpl},
and adding the bound $O(nD)$ on the number of incidences
between points and lines not assigned to the same component, we obtain the
bound asserted in the lemma.
\proofend

\medskip

\noindent{\bf Remark.}
Note that in both Lemma~\ref{incpl} and Lemma~\ref{incone}, we bound the number
of incidences between points on planar or conic components of $Z(f)$ and
\emph{all} the lines of $L_1$.

\paragraph{Pruning.}
To continue, we remove all the points of $P_1$ that lie in some planar or conic
component of $Z(f)$, and all the lines of $L_1$ that are fully contained in such
components. With the choice of $D=m^{1/2}/n^{1/4}$, we lose in the process
$$
O(m^{2/3}n^{1/3}s^{1/3}+m+nD) = O(m^{1/2}n^{3/4}+m^{2/3}n^{1/3}s^{1/3})
$$
incidences (recall that the term $m$ is subsumed by the term $m^{1/2}n^{3/4}$ for $m<n^{3/2}$).
Continue, for simplicity of notation, to denote the sets of remaining points and
lines as $P_1$ and $L_1$, respectively, and their sizes as $m$ and $n$.
Now each point is incident to at most $D^2$
lines (a fact that we will not use for this value of $D$), and no plane contains
more than $D$ lines of $L_1$, a crucial property for the next steps of the analysis.
That is, this allows us to replace the input parameter $s$, bounding the maximum
number of coplanar lines, by $D$; this is a key step that makes the induction work.

\paragraph{A new polynomial partitioning.}
We now return to the promised step of constructing a new polynomial partitioning.
We adapt the preceding notation, with a few modifications.
We choose a degree $E$, typically much smaller than $D$, and construct
a partitioning polynomial $g$ of degree $E$ for $P_1$. With an appropriate value
of $r=\Theta(E^3)$, we obtain $O(r)$ open cells, each containing at most $m/r$ points
of $P_1$, and each line of $L_1$ either crosses at most $E+1$ cells, or is fully
contained in $Z(g)$.

Set $P_2:= P_1\cap Z(g)$ and $P'_2:=P_1\setminus P_2$. Similarly, denote by
$L_2$ the set of lines of $L_1$ that are fully contained in $Z(g)$, and put
$L'_2:=L_1\setminus L_2$. We first dispose of incidences involving the lines
of $L_2$. (That is, now we first focus on incidences within $Z(g)$, and only then
turn to look at the cells.) By Lemma~\ref{incpl} and Lemma~\ref{incone}, the number
of incidences involving points $P_2$ that lie in some planar or conic component of
$Z(g)$, and all the lines of $L_2$, is
$$
O(m^{2/3}n^{1/3}s^{1/3}+m+nE) =
O(m^{1/2}n^{3/4}+m^{2/3}n^{1/3}s^{1/3}+n).
$$
(For $E\ll D$, this might be a gross overestimation, but we do not care.)
We remove these points from $P_2$, and remove all the lines of $L_2$ that are
contained in such components; continue to denote the sets of remaining points
and lines as $P_2$ and $L_2$. Now each point is incident to at most $E^2$ lines
of $L_2$ (Lemma~\ref{md2}), so the number of remaining incidences involving points
of $P_2$ is $O(mE^2)$; for $E$ suitably small, this bound will be subsumed by
$O(m^{1/2}n^{3/4})$.

Unlike the case of a ``large'' $D$, namely, $D=m^{1/2}/n^{1/4}$, here the
difficult part is to treat incidences within the cells of the partition.
Since $E\ll D$, we cannot use the naive bound $O(n^2+m)$ within each cell,
because that would make the overall bound too large. Therefore, to control the
incidence bound within the cells, we proceed in the following inductive manner.

For each cell $\tau$ of $\reals^3\setminus Z(g)$, put $P_\tau := P'_2\cap\tau$,
and let $L_\tau$ denote the set of the lines of $L'_2$ that cross $\tau$; put
$m_\tau = |P_\tau| \le m/r$, and $n_\tau = |L_\tau|$. Since every line
$\ell\in L_1$ (that is, of $L'_2$) crosses at
most $1+E$ components of $\reals^3\setminus Z(g)$, we have
$\sum_\tau n_\tau \le n(1+E)$.

It is important to note that at this point of the analysis the sizes of $P_1$
and of $L_1$ might be smaller than the original respective values $m$ and $n$.
In particular, we may no longer assume that $|P_1| > |L_1|^{\alpha_{j-1}}$, as
we did assume for $m$ and $n$. Nevertheless, in what follows $m$ and $n$ will
denote the original values, which serve as upper bounds for the respective
actual sizes of $P_1$ and $L_1$, and the induction will work correctly with
these values; see below for details.

In order to apply the induction hypothesis within the cells of the
partition, we want to assume that $m_\tau \le {n_\tau}^{\alpha_{j-1}}$ for
each $\tau$. To ensure that, we require that the number of lines of $L'_2$
that cross a cell be at most $n/E^2$. Cells $\tau$ that are crossed by
$\kappa n/E^2$ lines, for $\kappa>1$, are treated as if they occur
$\lceil \kappa \rceil$ times, where each incarnation involves all
the points of $P_\tau$, and at most $n/E^2$ lines of $L_\tau$. The number
of subproblems remains $O(E^3)$. Arguing similarly, we may also assume that
$m_\tau \le m/E^3$ for each cell $\tau$
(by ``duplicating'' each cell into a constant number of subproblems, if needed).

We therefore require that
${\displaystyle \frac {m} {E^3} \le \left(\frac{n}{E^2}\right)^{\alpha_{j-1}}}$.
(Note that, as already commented above, these are only upper bounds on the actual
sizes of these subsets, but this will have no real effect on the induction process.)
That is, we require
\begin{equation}
\label{Dlb}
E\ge \left( \frac {m} {n^{\alpha_{j-1}}} \right)^{1/(3-2\alpha_{j-1})} .
\end{equation}

With these preparations, we apply the induction hypothesis within each cell $\tau$,
recalling that no plane contains more than $D$ lines\footnote{%
  This was the main reason for carrying out the first partitioning step, as already noted.}
of $L'_2\subseteq L_1$, and get
\begin{align*}
I(P_\tau,L_\tau) & \le A_{j-1} \left( m_\tau^{1/2}n_\tau^{3/4} + m_\tau \right)
+ B\left( m_\tau^{2/3}n_\tau^{1/3}D^{1/3} + n_\tau \right) \\
& \le A_{j-1} \left( (m/E^3)^{1/2}(n/E^2)^{3/4} + m/E^3 \right)
+ B\left(+ (m/E^3)^{2/3}(n/E^2)^{1/3}D^{1/3} + n/E^2 \right) .
\end{align*}
Summing these bounds over the cells $\tau$, that is, multiplying them by
$O(E^3)$, we get, for a suitable absolute constant $b$,
$$
I(P'_2,L'_2) = \sum_{\tau} I(P_\tau,L_\tau)
\le b A_{j-1} \left( m^{1/2}n^{3/4} + m\right) + B\left(m^{2/3}n^{1/3}E^{1/3}D^{1/3} + nE \right) .
$$
We now require that $E=O(D)$. Then the last term satisfies
$nE = O(nD) = O(m^{1/2}n^{3/4})$, and, as already remarked, the preceding
term $m$ is also subsumed by the first term. The second term, after
substituting $D=O(m^{1/2}/n^{1/4})$, becomes
$O(m^{5/6}n^{1/4}E^{1/3})$.  Hence, with a slightly larger $b$, we have
$$
I(P'_2,L'_2) \le bA_{j-1} m^{1/2}n^{3/4} + bB m^{5/6}n^{1/4}E^{1/3} .
$$
Adding up all the bounds, including those for the portions of $P$ and $L$ that
were discarded during the first partitioning step, we obtain, for a suitable constant $c$,
$$
I(P,L) \le c\left(m^{1/2}n^{3/4} + m^{2/3}n^{1/3}s^{1/3} + n + mE^2\right)
+ bA_{j-1} m^{1/2}n^{3/4} + bB m^{5/6}n^{1/4}E^{1/3} .
$$

We choose $E$ to ensure that the two $E$-dependent terms are dominated by
the term $m^{1/2}n^{3/4}$. That is,
\begin{align*}
m^{5/6}n^{1/4}E^{1/3} & \le m^{1/2}n^{3/4} ,\quad\text{or}\quad E\le n^{3/2}/m , \\
\text{and} \quad mE^2 & \le m^{1/2}n^{3/4} ,\quad\text{or}\quad E\le n^{3/8}/m^{1/4} .
\end{align*}
Since $n^{3/2}/m = \left( n^{3/8}/m^{1/4} \right)^4$, and both sides
are $\ge 1$, the latter condition is stricter, and we ignore the former.
As already noted, we also require that $E=O(D)$; specifically, we require
that $E\le m^{1/2}/n^{1/4}$.

In conclusion, recalling~(\ref{Dlb}), the two constraints on the choice of $E$ are
\begin{equation} \label{econ}
\left( \frac {m} {n^{\alpha_{j-1}}} \right)^{1/(3-2\alpha_{j-1})} \le E
\le \min\left\{ \frac{n^{3/8}}{m^{1/4}} , \frac{m^{1/2}}{n^{1/4}} \right\} ,
\end{equation}
and, for these constraints to be compatible, we require that
$$
\left( \frac {m} {n^{\alpha_{j-1}}} \right)^{1/(3-2\alpha_{j-1})}
\le \frac{n^{3/8}}{m^{1/4}} ,\quad\quad\text{or}\quad\quad
m \le n^{\frac{9+2\alpha_{j-1}}{2(7-2\alpha_{j-1})}} ,
$$
and that
$$
\left( \frac {m} {n^{\alpha_{j-1}}} \right)^{1/(3-2\alpha_{j-1})}
\le \frac{m^{1/2}}{n^{1/4}} ,
$$
which fortunately always holds, as is easil;y checked, since $m\le n^{3/2}$
and $\alpha_{j-1} \ge 1/2$. Note that we have not explicitly stated any
concrete choice of $E$; any value satisfying (\ref{econ}) will do. We put
$$
\alpha_j := \frac{9+2\alpha_{j-1}}{2(7-2\alpha_{j-1})} ,
$$
and conclude that if $m\le n^{\alpha_j}$ then the bound asserted in the theorem
holds, with $A_j = bA_{j-1} + c$ and $B=c$. This completes the induction step.
Note that the recurrence $A_j = bA_{j-1}+c$ solves to $A_j = O(b^j)$.

It remains to argue that the induction covers the entire range $m = O(n^{3/2})$.
Using the above recurrence for the $\alpha_j$'s, with $\alpha_0 = 1/2$,
it easily follows that
$$
\alpha_j = \frac32 - \frac{2}{j+2} ,
$$
for each $j\ge 0$, showing that $\alpha_j$ converges to $3/2$, implying
that the entire range $m = O(n^{3/2})$ is covered by the induction.

To calibrate the dependence of the constant of proportionality on $m$ and $n$,
we note that, for $n^{\alpha_{j-1}}\le m < n^{\alpha_j}$, the constant is $O(b^j)$. We have
$$
\frac32 - \frac{2}{j+1} = \alpha_{j-1} \le \frac{\log m}{\log n} ,
\quad\quad\text{or}\quad\quad
j \le \frac{\frac12 + \frac{\log m}{\log n}}{\frac32 - \frac{\log m}{\log n}} =
\frac{\log (m^2n)}{\log (n^3/m^2)} .
$$
This establishes the expression for $A_{m,n}$ given in the statement of the theorem.

\paragraph{Handling the middle ground $m\approx n^{3/2}$.}
Some care is needed when $m$ approaches $n^{3/2}$, because of the
potentially unbounded growth of the constant $A_j$. To handle this situation,
we simply fix a value $j$, in the manner detailed below, write $m=kn^{\alpha_j}$,
solve $k$ separate problems, each involving $m/k=n^{\alpha_j}$ points of $P$
and all the $n$ lines of $L$, and sum up the resulting incidence bounds. We then get
\begin{align*}
I(P,L) & \le akb^j \left( (m/k)^{1/2}n^{3/4} + (m/k) \right)
+ kB \left( (m/k)^{2/3}n^{1/3}s^{1/3} + n \right ) \\
& = ak^{1/2}b^j m^{1/2}n^{3/4} + ab^j m + k^{1/3}B m^{2/3}n^{1/3}s^{1/3} + kB n ,
\end{align*}
for a suitable absolute constant $a$.
Recalling that $\alpha_j = \frac32 - \frac{2}{j+2}$, we have
$$
k \le m/n^{\alpha_j} \le n^{3/2}/n^{\alpha_j} = n^{2/(j+2)} .
$$
Hence the coefficient of the leading term in the above bound is bounded by
$an^{1/(j+2)}b^j$, and we (asymptotically) minimize this expression by
choosing
$$
j = j_0 := \sqrt{\log n}/\sqrt{\log b} .
$$
With this choice all the other
coefficients are also dominated by the leading coefficient, and we obtain
\begin{equation} \label{eq:mn32}
I(P,L) = O\left( 2^{2\sqrt{\log b}\sqrt{\log n}} \left(
m^{1/2}n^{3/4} + m^{2/3}n^{1/3}s^{1/3} + m + n \right) \right) .
\end{equation}
In other words, the bound in (\ref{st}) and (\ref{amn}) holds for any $m\le n^{3/2}$,
but, for $m \ge n^{\alpha_{j_0}}$ one should use instead the bound in (\ref{eq:mn32}),
which controls the exponential growth of the constants of proportionality within this range.

\paragraph{The case $m > n^{3/2}$.}
The analysis of this case is, in a sense, a mirror image of the preceding analysis,
except for a new key lemma (Lemma~\ref{nd2}).
For the sake of completeness, we repeat a sizeable portion of the analysis, providing
many of the relevant (often differing) details.

We partition this range into a sequence of ranges $m\ge n^{\alpha_0}$,
$n^{\alpha_1} \le m < n^{\alpha_0},\ldots$, where $\alpha_0 = 2$
and the sequence $\{\alpha_j\}_{j\ge 0}$ is decreasing and converges
to $3/2$. The induction is on the index $j$ of the range
$n^{\alpha_{j}} \le m < n^{\alpha_{j-1}}$, and establishes (\ref{st})
for $m$ in this range, with a coefficient $A_j$ (written in
(\ref{st},\ref{amn}) as $A_{m,n}$) that increases with $j$.

The base range of the induction is $m\ge n^2$, where the
trivial general upper bound on point-line incidences in any
dimension, dual to the one used in the previous case, yields
$I=O(n^2+m)=O(m)$, so (\ref{st}) holds for a
sufficiently large choice of the initial constant $A_0$.

Assume then that (\ref{st}) holds for all $m \ge n^{\alpha_{j-1}}$
for some $j \ge 1$, and consider an instance of the problem with
$n^{3/2} \le m < n^{\alpha_{j-1}}$ (again, the lower bound will increase,
to $n^{\alpha_j}$, to facilitate the induction step).

For a parameter $r$, to be specified later,
apply the polynomial partition theorem to obtain an $r$-partitioning trivariate
(real) polynomial $f$ of degree $D=O(r^{1/3})$. That is, every connected
component of $\reals^3\setminus Z(f)$ contains at most $m/r$ points of $P$, and
the number of components of $\reals^3\setminus Z(f)$ is $O(D^3) = O(r)$.

Set $P_1:= P\cap Z(f)$ and $P'_1:=P\setminus P_1$. Each line $\ell\in L$ is
either fully contained in $Z(f)$ or intersects it in at most $D$ points.
Let $L_1$ denote the subset of lines of $L$ that are
fully contained in $Z(f)$ and put $L'_1 = L\setminus L_1$.
As before, we have
$$
I(P,L) = I(P_1,L_1) + I(P_1,L'_1) + I(P'_1,L'_1) .
$$
We have
$$
I(P_1,L'_1) \le |L'_1|\cdot D \le nD ,
$$
and we estimate $I(P'_1,L'_1)$ as follows. For each cell $\tau$ of
$\reals^3\setminus Z(f)$, put $P_\tau = P\cap \tau$ (that is,
$P'_1\cap\tau$), and let $L_\tau$ denote the set of the lines
of $L'_1$ that cross $\tau$; put $m_\tau = |P_\tau| \le m/r$, and
$n_\tau = |L_\tau|$. As before, we have
$\sum_\tau n_\tau \le n(1+D)$, so the average number of lines
that cross a cell is $O(n/D^2)$. Arguing as above,
we may assume, by possibly increasing the number of cells by a constant factor, that each
$n_\tau$ is at most $n/D^2$. Clearly, we have
$$
I(P'_1,L'_1) = \sum_{\tau} I(P_\tau,L_\tau).
$$
For each $\tau$ we use the trivial dual bound, mentioned above,
$I(P_\tau,L_\tau) = O(n_\tau^2+m_\tau)$. Summing over the cells, we get
$$
I(P'_1,L'_1) = \sum_{\tau} I(P_\tau,L_\tau) =
O\left(D^3\cdot(n/D^2)^2 + m \right) = O\left(n^2/D + m \right) .
$$
For the initial value of $D$, we take $D = n^2/m$, noting that
$1\le D^3\le m$ because $n^{3/2}\le m\le n^2$, and get the bound
$$
I(P'_1,L'_1) + I(P_1,L'_1) = O(n^2/D+m+nD) = O(m + n^3/m) = O(m) ,
$$
where the latter bound follows since $m\ge n^{3/2}$.

It remains to estimate $I(P_1,L_1)$. Since all the incidences involving
any point in $P'_1$ and/or any line in $L'_1$ have been accounted for, we
discard these sets, and remain with $P_1$ and $L_1$ only. As before, we forget
the preceding polynomial partitioning step, and start afresh, applying a
new polynomial partitioning to $P_1$ with a polynomial $g$ of degree $E$,
which will typically be much smaller than $D$, but still non-constant.

For this case we need the following lemma, which can be regarded, in some sense,
as a dual (albeit somewhat more involved) version of Lemma~\ref{md2}.
Unlike the rest of the analysis, the best way to prove this lemma is by
switching to the complex projective setting. This is needed for one key step
in the proof, where we need the property that the projection of a complex projective
variety is a variety. Once this is done, we can switch back to the real affine case,
and complete the proof.

Here is a very quick review of the transition to the complex projective setup.
A real affine algebraic variety $X$, defined by a collection of real polynomials,
can also be regarded as a complex projective variety. (Technically, one needs
to take the \emph{projective closure} of the \emph{complexification} of $X$;
details about these standard operations can be found, e.g., in
Bochnak et al.~\cite[Proposition 8.7.17]{BCR}
and in Cox et al.~\cite[Definition 8.4.6]{CLO}.) If $f$ is an irreducible
polynomial over $\reals$, it might still be reducible over $\cplx$, but then
it must have the form $f=g\bar g$,
where $g$ is an irreducible complex polynomial and $\bar g$ is its complex
conjugate. (Indeed, if $h$ is any irreducible factor of $f$, then $\bar h$
is also an irreducible factor of $f$, and therefore $h\bar h$ is a real
polynomial dividing $f$. As $f$ is irreducible over $\reals$, the claim follows.)

In the following lemma, adapting a notation used in earlier works,
we say that a point $p\in P_1$ is \emph{$1$-poor} (resp.,
\emph{$2$-rich}) if it is incident to at most one line (resp., to at
least two lines) of $L_1$.

Recall also that a \emph{regulus} is a doubly-ruled surface in
$\reals^3$ or in $\cplx^3$. It is the union of all lines that pass through
three fixed pairwise skew lines; it is a quadric, which is either a
hyperbolic paraboloid or a one-sheeted hyperboloid.

\begin{lemma} \label{nd2}
Let $f$ be an irreducible polynomial in $\cplx[x,y,z]$, such that $Z(f)$ is not
a complex plane nor a complex regulus, and let $L_1$ be a finite set of lines 
fully contained in $Z(f)$. Then, with the possible exception of at most two lines,
each line $\ell\in L_1$ is incident to at most $O(D^3)$ $2$-rich points.
\end{lemma}
\noindent{\bf Proof.}
The strategy of the proof is to charge each incidence of $\ell$ with some
$2$-rich point $p$ to an intersection of $\ell$ with another line
of $L_1$ that passes through $p$, and to argue that, in general, there can
be only $O(D^3)$ such other lines. This in turn will be shown by arguing that
the union of all the lines that are fully contained in $Z(f)$ and pass through
$\ell$ is a one-dimensional variety, of degree $O(D^3)$, from which the claim
will follow. As we will show, this will indeed be the case except when $\ell$
is one of at most two ``exceptional'' lines on $Z(f)$.

Fix a line $\ell$ as in the lemma, assume for simplicity that it passes through
the origin, and write it as $\{tv_0 \mid t\in\cplx\}$; since $\ell$ is a real
line, $v_0$ can be assumed to be real. Consider the union $V(\ell)$ of all the
lines that are fully contained in $Z(f)$ and are incident to $\ell$; that is,
$V(\ell)$ is the union of $\ell$ with the set of all points $p\in Z(f)\setminus\ell$
for which there exists $t\in\cplx$ such that the line connecting $p$ to
$tv_0\in\ell$ is fully contained in $Z(f)$. In other
words, for such a $t$ and for each $s\in\cplx$, we have $f((1-s)p+stv_0)=0$.
Regarding the left-hand side as a polynomial in $s$, we can write it as
${\displaystyle \sum_{i=0}^D G_i(p;t) s^i \equiv 0}$, for suitable (complex)
polynomials $G_i(p;t)$ in $p$ and $t$, each of total degree at most $D$.
In other words, $p$ and $t$ have to satisfy the system
\begin{equation} \label{fsyst}
G_0(p;t) = G_1(p;t) = \cdots = G_D(p;t) = 0 ,
\end{equation}
which defines an algebraic variety $\sigma(\ell)$ in $\P^4(\cplx)$.
Note that, substituting $s=0$, we have $G_0(p;t)\equiv f(p)$, and that the
limit points $(tv_0,t)$ (corresponding to points on $\ell$) also satisfy
this system, since in this case $f((1-s)tv_0+stv_0)=f(tv_0)=0$ for all $s$.

In other words, $V(\ell)$ is the projection of $\sigma(\ell)$ into $\P^3(\cplx)$,
given by $(p,t)\mapsto p$. For each $p\in Z(f)\setminus\ell$ this system has only
finitely many solutions in $t$, for otherwise the plane spanned by $p$ and $\ell_0$
would be fully contained in $Z(f)$, contrary to our assumption.

By the projective extension theorem (see, e.g., \cite[Theorem 8.6]{CLO}),
the projection of $\sigma(\ell)$ into $\P^3(\cplx)$, in which $t$ is discarded,
is an algebraic variety $\tau(\ell)$. We observe that $\tau(\ell)$ is contained
in $Z(f)$, and is therefore of dimension at most two.

Assume first that $\tau(\ell)$ is two-dimensional. As $f$ is irreducible
over $\cplx$, we must have $\tau(\ell) = Z(f)$.  This implies that each
point $p\in Z(f)\setminus\ell$ is incident to a (complex) line that is fully contained
in $Z(f)$ and is incident to $\ell$. In particular, $Z(f)$ is ruled by
complex lines.

By assumption, $Z(f)$ is neither a complex plane nor a complex
regulus. We may also assume that $Z(f)$ is not a complex cone, for
then each line in $L_1$ is incident to at most one 2-rich point
(namely, the apex of $Z(f)$), making the assertion of the lemma
trivial. It then follows that $Z(f)$ is an irreducible singly ruled
(complex) surface.
As argued in Guth and Katz~\cite{GK2} (see also our
companion paper~\cite{SS3dvar} for an independent analysis of this
situation, which caters more explicitly to the complex setting too),
$Z(f)$ can contain at most two lines $\ell$ with this property.

Excluding these (at most) two exceptional lines $\ell$, we may thus assume
that $\tau(\ell)$ is (at most) a one-dimensional curve.

Clearly, by definition, each point $(p,t)\in\sigma(\ell)$, except
for $p\in\ell$, defines a line $\lambda$, in the original 3-space,
that connects $p$ to $tv_0$, and each point $q\in\lambda$ satisfies
$(q,t)\in\sigma(\ell)$. Hence, the line $\{(q,t) \mid q\in\lambda\}$
is fully contained in $\sigma(\ell)$, and therefore the line
$\lambda$ is fully contained in $\tau(\ell)$. Since $\tau(\ell)$ is
one-dimensional, this in turn implies
(see, e.g., \cite[Lemma 2.3]{SS4d})
that $\tau(\ell)$ is a \emph{finite} union of (complex)
lines, whose number is at most $\deg (\tau(\ell))$. This also
implies that $\sigma(\ell)$ is the union of the same number of
lines, and in particular $\sigma(\ell)$ is also one-dimensional,
and the number of lines that it contains is at most $\deg (\sigma(\ell))$.

We claim that this latter degree is at most $O(D^3)$. This follows
from a well-known result in algebra (see, e.g., Schmid \cite[Lemma
2.2]{Schmid}), that asserts that, since $\sigma(\ell)$ is a
one-dimensional curve in $\P^4(\cplx)$, and is the common zero set
of polynomials, each of degree $O(D)$, its degree is $O(D^3)$.

This completes the proof of the lemma. (The passage from the
complex projective setting back to the real affine one is trivial
for this property.) \proofend

\begin{corollary} \label{nd3}
Let $f$ be a real or complex trivariate polynomial of degree $D$, such that
(the complexification of) $Z(f)$ does not contain any complex plane nor any
complex regulus. Let $L_1$ be a set of $n$ lines fully contained in $Z(f)$, 
and let $P_1$ be a set of $m$ points contained in $Z(f)$. Then
$I(P_1,L_1) = O(m+nD^3)$.
\end{corollary}
\noindent{\bf Proof.}
Write $f=\prod_{i=1}^s f_i$ for its decomposition into irreducible factors, 
for $s\le D$. We apply Lemma~\ref{nd2} to each complex factor $f_i$ of the 
$f$. By the observation preceding Lemma~\ref{nd2},some of these factors might 
be complex (non-real) polynomials, even when $f$ is real. That is, regardless 
of whether the original $f$ is real or not, we carry out the analysis in
the complex projective space $\P^3(\cplx)$, and regard $Z(f_i)$ as
a variety in that space.

Note also that, by focussing on the single irreducible component $Z(f_i)$ of 
$Z(f)$, we consider only points and lines that are fully contained in $Z(f_i)$.
We thus shrink $P_1$ and $L_1$ accordingly, and note that the notions of being 
2-rich or 1-poor are now redefined with respect to the reduced sets. 
All of this will be rectified at the end of the proof.

Assign each line $\ell\in L_1$ to the first component $Z(f_i)$, in the above order,
that fully contains $\ell$, and assign each point $p\in P_1$ to the
first component that contains it. If a point $p$ and a line $\ell$ are incident,
then either they are both assigned to the same component $Z(f_i)$, or $p$ is
assigned to some component $Z(f_i)$ and $\ell$, which is assigned to a later
component, is not contained in $Z(f_i)$. Each incidence of the latter kind
can be charged to a crossing between $\ell$ and $Z(f_i)$, and the total
number of these crossings is $O(nD)$. It therefore suffices to consider
incidences between points and lines assigned to the same component.
Moreover, if a point $p$ is 2-rich with respect to the entire collection $L_1$
but is 1-poor with respect to the lines assigned to its component, then all
of its incidences except one are accounted by the preceding term $O(nD)$,
which thus takes care also of the single incidence within $Z(f_i)$.

By Lemma~\ref{nd2}, for each $f_i$, excluding at most two exceptional lines,
the number of incidences between a line assigned to (and contained in) $Z(f_i)$ 
and the points assigned to $Z(f_i)$ that are still 2-rich within $Z(f_i)$, is
$O(\deg(f_i)^3) = O(D^3)$.  Summing over all relevant lines, we get the bound $O(nD^3)$.

Finally, each irreducible component $Z(f_i)$ can contain at most two exceptional 
lines, for a total of at most $2D$ such lines. The number of $2$-rich points on each
such line $\ell$ is at most $n$, since each such point is incident to another line, 
so the total number of corresponding incidences is at most $O(nD)$, which is subsumed
by the preceding bound $O(nD^3)$. The number of incidences with $1$-poor points is,
trivially, at most $m$. This completes the proof of the corollary. \proofend

\paragraph{Pruning.}
In the preceding lemma and corollary, we have excluded planar and reguli components of $Z(f)$.
Arguing as in the case of small $m$, the number of incidences involving points that
lie on planar components of $Z(f)$ is $O(m^{2/3}n^{1/3}s^{1/3}+m)$ (see Lemma~\ref{incpl}),
and the number of incidences involving points that lie on conic components of $Z(f)$ is
$O(m+nD) = O(m)$ (see Lemma~\ref{incone}).
A similar bound holds for points on the reguli components. Specifically, we assign
each point and line to a regulus that contain them, if one exists, in the same
first-come first-serve manner used above. Any point $p$ can be incident to at most
two lines that are fully contained in the regulus to which it is assigned, and
any other incidence of $p$ with a line $\ell$ can be uniquely charged to the intersection
of $\ell$ with that regulus, for a total (over all lines and reguli) of $O(nD)$ incidences.

We remove all points that lie in any such component and all lines that are fully contained
in any such component.  With the choice of $D=n^2/m$, we lose in the process
$$
O(m^{2/3}n^{1/3}s^{1/3}+m+nD) = O(m + m^{2/3}n^{1/3}s^{1/3})
$$
incidences (recall that $nD\le m$ for $m\ge n^{3/2}$). For the remainder sets, which we
continue to denote as $P_1$ and $L_1$, respectively, no plane contains
more than $O(D)$ lines of $L_1$, as argued in Lemma~\ref{donpl}.

\paragraph{A new polynomial partitioning.}
We adapt the notation used in the preceding case, with a few modifications.
We choose a degree $E$, typically much smaller than $D$, and construct
a partitioning polynomial $g$ of degree $E$ for $P_1$. With an appropriate value
of $r=\Theta(E^3)$, we obtain $O(r)$ cells, each containing at most $m/r$ points
of $P_1$, and each line of $L_1$ either crosses at most $E+1$ cells, or is fully
contained in $Z(g)$.

Set $P_2:= P_1\cap Z(g)$ and $P'_2:=P_1\setminus P_2$. Similarly, denote by
$L_2$ the set of lines of $L_1$ that are fully contained in $Z(g)$, and put
$L'_2:=L_1\setminus L_2$. We first dispose of incidences involving the lines
of $L_2$. By Lemma~\ref{incpl} and the preceding arguments, the number of
incidences involving points of $P_2$ that lie in some planar, conic, or regulus
component of $Z(g)$, and all the lines of $L_2$, is
$$
O(m^{2/3}n^{1/3}s^{1/3}+m+nE) .
$$
We remove these points from $P_2$, and remove all the lines of $L_2$ that are
contained in such components. Continue to denote the sets of remaining points
and lines as $P_2$ and $L_2$. By Corollary~\ref{nd3}, the number of incidences
between $P_2$ and $L_2$ is $O(m+nE^3)$.

To complete the estimation, we need to bound the number of incidences
in the cells of the partition, which we do inductively, as before. Specifically,
for each cell $\tau$ of $\reals^3\setminus Z(g)$, put $P_\tau := P'_2\cap\tau$,
and let $L_\tau$ denote the set of the lines of $L'_2$ that cross $\tau$; put
$m_\tau = |P_\tau| \le m/r$, and $n_\tau = |L_\tau|$. Since every line $\ell\in L_0$ crosses at
most $1+E$ components of $\reals^3\setminus Z(g)$, we have
$\sum_\tau n_\tau \le n(1+E)$, and, arguing as above, we may assume that each
$n_\tau$ is at most $n/E^2$, and each $m_\tau$ is at most $m/E^3$.
To apply the induction hypothesis in each cell, we therefore require that
${\displaystyle \frac {m} {E^3} \ge \left(\frac{n}{E^2}\right)^{\alpha_{j-1}}}$.
(As before, the actual sizes of $P_1$ and $L_1$ might be smaller than the respective
original values $m$ and $n$. We use here the original values, and note, similar
to the preceding case, that the fact that these are only upper bounds on the
actual sizes is harmless for the induction process.) That is, we require
\begin{equation}
\label{Dlb2}
E\ge \left( \frac {n^{\alpha_{j-1}}}{m} \right)^{1/(2\alpha_{j-1}-3)} .
\end{equation}

With these preparations, we apply the induction hypothesis within each cell $\tau$,
recalling that no plane contains more than $D$ lines of $L'_2\subseteq L_1$, and get
\begin{align*}
I(P_\tau,L_\tau) & \le A_{j-1} \left( m_\tau^{1/2}n_\tau^{3/4}+ m_\tau \right)
+ B\left( m_\tau^{2/3}n_\tau^{1/3}D^{1/3} + n_\tau \right) \\
& \le A_{j-1} \left( (m/E^3)^{1/2}(n/E^2)^{3/4}+ m/E^3 \right)
+ B\left( (m/E^3)^{2/3}(n/E^2)^{1/3}D^{1/3} + n/E^2 \right) .
\end{align*}
Summing these bounds over the cells $\tau$, that is, multiplying them by
$O(E^3)$, we get, for a suitable absolute constant $b$,
$$
I(P'_2,L'_2) = \sum_{\tau} I(P_\tau,L_\tau)
\le b A_{j-1} \left( m^{1/2}n^{3/4} + m\right) +
bB\left(m^{2/3}n^{1/3}E^{1/3}D^{1/3} + nE \right) .
$$
Requiring that $E\le m/n$, the last term satisfies $nE\le m$,
and the first term is also at most $O(m)$ (because $m\ge n^{3/2}$).
The second term, after substituting $D=O(n^2/m)$, becomes
$O(m^{1/3}nE^{1/3})$.  Hence, with a slightly larger $b$, we have
$$
I(P'_2,L'_2) \le bA_{j-1} m + bB m^{1/3}nE^{1/3} .
$$
Collecting all partial bounds obtained so far, we obtain
$$
I(P,L) \le c\left(m^{2/3}n^{1/3}s^{1/3} + m + nE^3\right)
+ bA_{j-1} m + bB m^{1/3}nE^{1/3} ,
$$
for a suitable constant $c$.
We choose $E$ to ensure that the two $E$-dependent terms are dominated by $m$.
That is,
$$
m^{1/3}nE^{1/3} \le m ,\quad\text{or}\quad E\le m^{2}/n^{3} , \quad\quad\text{and}\quad\quad
nE^3 \le m ,\quad\text{or}\quad E\le m^{1/3}/n^{1/3} .
$$
In addition, we also require that $E\le m/n$, but, as is easily seen, both
of the above constraints imply that $E\le m/n$, so we get this latter
constraint for free, and ignore it in what follows.

As is easily checked, the second constraint
$E\le m^{1/3}/n^{1/3}$ is stricter than the first constraint
$E\le m^2/n^3$ for $m\ge n^{8/5}$, and the situation is reversed when
$m\le n^{8/5}$. So in our inductive descent of $m$, we first consider
the second constraint, and then switch to the first constraint.

Hence, in the first part of this analysis,
the two constraints on the choice of $E$ are
$$
\left( \frac {n^{\alpha_{j-1}}}{m} \right)^{1/(2\alpha_{j-1}-3)} \le E
\le \frac{m^{1/3}}{n^{1/3}} ,
$$
and, for these constraints to be compatible, we require that
$$
\left( \frac {n^{\alpha_{j-1}}}{m} \right)^{1/(2\alpha_{j-1}-3)}
\le \frac{m^{1/3}}{n^{1/3}} , \quad\quad\text{or}\quad\quad
m \ge n^{\frac{5\alpha_{j-1}-3}{2\alpha_{j-1}}} .
$$
We start the process with $\alpha_0=2$, and take
${\displaystyle \alpha_1 := \frac{5\alpha_{0}-3}{2\alpha_{0}} = 7/4}$.
As this is still larger than $8/5$, we perform two additional rounds of the induction,
using the same constraints, leading to the exponents
$$
\alpha_2 = \frac{5\alpha_1-3}{2\alpha_1} = \frac{23}{14}, \quad\text{and}\quad
\alpha_3 = \frac{5\alpha_2-3}{2\alpha_2} = \frac{73}{46} < \frac{8}{5} .
$$
To play it safe, we reset $\alpha_3:=8/5$, and establish the induction
step for $m\ge n^{8/5}$. We can then proceed to the second part, where
the two constraints on the choice of $E$ are
$$
\left( \frac {n^{\alpha_{j-1}}}{m} \right)^{1/(2\alpha_{j-1}-3)} \le E
\le \frac{m^{2}}{n^{3}} ,
$$
and, for these constraints to be compatible, we require that
$$
\left( \frac {n^{\alpha_{j-1}}}{m} \right)^{1/(2\alpha_{j-1}-3)}
\le \frac{m^{2}}{n^{3}} , \quad\quad\text{or}\quad\quad
m \ge n^{\frac{7\alpha_{j-1}-9}{4\alpha_{j-1}-5}} .
$$
We define, for $j\ge 4$,
${\displaystyle \alpha_j = \frac{7\alpha_{j-1}-9}{4\alpha_{j-1}-5}}$.
Substituting $\alpha_3=8/5$ we get $\alpha_4=11/7$, and in general
a simple calculation shows that
$$
\alpha_j = \frac32 + \frac{1}{4j-2} ,
$$
for $j\ge 3$. This sequence does indeed converge to $3/2$ as $j\to\infty$,
implying that the entire range $m = \Omega(n^{3/2})$ is covered by the induction.

In both parts, we conclude that if $m\ge n^{\alpha_j}$ then the bound asserted in
the theorem holds with $A_j = bA_{j-1}+c$,, and $B=c$. This completes the induction step.

Finally, we calibrate the dependence of the constant of proportionality on $m$ and $n$,
by noting that, for $n^{\alpha_{j}}\le m < n^{\alpha_{j-1}}$, the constant is $O(b^j)$. We have
$$
\frac32 + \frac{1}{4j-6} = \alpha_{j-1} \ge \frac{\log m}{\log n} , \quad\text{or}\quad
j \le \frac{3 \frac{\log m}{\log n}-4}{2\frac{\log m}{\log n}-3} =
\frac{\log\left( m^3/n^4\right)}{\log\left(m^2/n^3\right)} .
$$
(Technically, this only handles the range $j\ge 3$, but, for an asymptotic bound, we can
extend it to $j=1,2$ too.)
This establishes the explicit expression for $A_{m,n}$ for this range, as stated
in the theorem, and completes its proof.
\proofend

Again, as in the case of a small $m$, we need to be careful when $m$ approaches
$n^{3/2}$. Here we can fix a $j$, assume that $n^{3/2} \le m < n^{\alpha_j}$,
and set $k:= m/n^{\alpha'_j}$, where $\alpha'_j = 3/2 - 2/(j+2)$ is the $j$-th
index in the hierarchy for $m\le n^{3/2}$. That is,
$$
k\le n^{\alpha_j-\alpha'_j} = \frac{1}{4j-2} + \frac{2}{j+2} .
$$
As before, we now solve $k$ separate
subproblems, each with $m/k$ points of $P$ and all the lines of $L$, and sum up
the resulting incidence bounds. The analysis is similar to the one used above,
and we omit its details. It yields almost the same bound as in (\ref{eq:mn32}),
where the slightly larger upper bound on $k$ leads to the slightly larger bound
$$
I(P,L) = O\left( 2^{\sqrt{4.5}\sqrt{\log b}\sqrt{\log n}} \left(
m^{1/2}n^{3/4} + m^{2/3}n^{1/3}s^{1/3} + m + n \right) \right) ,
$$
with a slightly different absolute constant $b$.


\section{Discussion}

In this paper we derived an asymptotically tight bound for the
number of incidences between a set $P$ of points and a set $L$ of
lines in $\reals^3$. This bound has already been established by Guth and
Katz~\cite{GK2}, where the main tool was the use of partitioning
polynomials. As already mentioned, the main novelty here is to use two 
separate partitioning polynomials of different degrees; the one with the higher degree
is used as a pruning mechanism, after which the maximum number of coplanar
lines of $L$ can be better controlled (by the degree $D$ of the polynomial),
which is a key ingredient in making the inductive argument work.

The second main tool of Guth and Katz was the Cayley--Salmon
theorem. This theorem says that a surface in $\reals^3$ of degree
$\Deg$ cannot contain more than $11\Deg^2-24\Deg$ lines, unless it
is \emph{ruled by lines}. This is an ``ancient'' theorem, from the
19th century, combining algebraic and differential geometry, and its
re-emergenece in recent years has kindled the interest of the
combinatorial geometry community in classical (and modern) algebraic
geometry. New proofs of the theorem were obtained (see, e.g., Terry
Tao's blog~\cite{tt:blog}), and generalizations
to higher dimensions have also been developed (see
Landsberg~\cite{Land}). However, the theorem only holds over the
complex field, and using it over the reals requires some care.

There is also an alternative way to bound the number of point-line incidences
using flat and singular points. However, as already remarked, these two,
as well as the Cayley--Salmon machinery, are non-trivial constructs,
especially in higher dimensions, and their generalization to other problems
in combinatorial geometry (even incidence problems with curves other
than lines or incidences with lines in higher dimensions) seem quite
difficult (and are mostly open). It is therefore
of considerable interest to develop alternative, more elementary
interfaces between algebraic and combinatorial geometry, which is a
primary goal of the present paper (as well as of Guth's
recent work~\cite{Gu14}). 

In this regard, one could perhaps view Lemma~\ref{md2} and 
Corollary~\ref{nd3} as certain weaker analogs of the Cayley--Salmon 
theorem, which are nevertheless easier to derive, without having to
use differential geometry. Some of the tools in Guth's paper~\cite{Gu14}
might also be interpreted as such weaker variants of the Cayley--Salmon
theory. It would be interesting to see suitable extensions of these tools
to higher dimensions.

Besides the intrinsic interest in simplifying the Guth--Katz analysis,
the present work has been motivated by our study of incidences between
points and lines in four dimensions. This has begun in a year-old
companion paper~\cite{surf-socg}, where we have used the
the polynomial partitioning method, with a polynomial of constant
degree. This, similarly to Guth's work in three dimensions~\cite{Gu14},
has resulted in a slightly weaker bound and considerably stricter
assumptions concerning the input set of lines. In a more involved follow-up
study~\cite{SS4d}, we have managed to improve the bound, and to get rid
of the restrictive assumptions, using two partitioning steps, with
polynomials of non-constant degrees, as in the present paper.
However, the analysis in \cite{SS4d} is not as simple as in the present paper,
because, even though there are generalizations of the Cayley--Salmon
theorem to higher dimensions (due to Landsberg, as mentioned above),
it turns out that a thorough investigation of the variety of lines
fully contained in a given hypersurface of non-constant degree, is a
fairly intricate and challenging problem, raising many deep questions
in algebraic geometry, some of which are still unresolved.

One potential application of the techniques used in this paper, mainly
the interplay between partitioning polynomials of different degrees,
is to the problem, recently studied by
Sharir, Sheffer and Zahl~\cite{SSZ}, of bounding the number of
incidences between points and circles in $\reals^3$. That paper
uses a partitioning polynomial of constant degree, and,
as a result, the term that caters to incidences
within lower-dimensional spaces (such as our term $m^{2/3}n^{1/3}s^{1/3}$)
does not go well through the induction mechanism, and consequently the bound
derived in \cite{SSZ} was weaker. We believe that our technique can improve
the bound of \cite{SSZ} in terms of this ``lower-dimensional'' term.

A substantial part of the present paper (half of the proof of the theorem)
was devoted to the treatment of the case $m>n^{3/2}$. However, under the
appropriate assumptions, the number of points incident to at least two
lines was shown by Guth and Katz~\cite{GK2} to be bounded by $O(n^{3/2})$.
A recent note by Koll\'ar~\cite{Kol} gives a simplified proof, including
an explicit multiplicative constant. In his work, Koll\'ar does not use
partitioning polynomials, but employs more advanced algebraic geometric
tools, like the \emph{arithmetic genus} of a curve, which serves as an upper
bound for the number of singular points. If we accept (pedagogically) the upper
bound $O(n^{3/2})$ for the number of 2-rich points as a ``black box'',
the regime in which $m>n^{3/2}$ becomes irrelevant, and can be discarded
from the analysis, thus greatly simplifying the paper.

A challenging problem is thus to find an elementary proof that the number
of points incident to at least two lines is $O(n^{3/2})$ (e.g., without
the use of the Cayley--Salmon theorem or the tools used by Koll\'ar).
Another challenging (and probably harder) problem is to improve the bound
of Guth and Katz when the bound $s$ on the maximum number of mutually
coplanar lines is $\ll n^{1/2}$: In their original derivation, Guth and Katz~\cite{GK2}
consider mainly the case $s=n^{1/2}$, and the lower bound constrcution in
\cite{GK2} also has $s=n^{1/2}$.
Another natural further research direction is to find further applications
of partitioning polynomials of intermediate degrees.


\begin{thebibliography}{99}


\bibitem{BS}
S. Basu and M. Sombra,
Polynomial partitioning on varieties and point-hypersurface incidences in four dimensions,
in arXiv:1406.2144.

\bibitem{BCR}
J. Bochnak, M. Coste and M. F. Roy,
{\it Real Algebraic Geometry},
Springer Verlag, Heidelberg, 1998.

\bibitem{CEGSW}
K. Clarkson, H. Edelsbrunner, L. Guibas, M. Sharir and E. Welzl,
Combinatorial complexity bounds for arrangements of curves and spheres,
{\it Discrete Comput. Geom.} 5 (1990), 99--160.

\bibitem{CLO}
D. Cox, J. Little and D. O'Shea,
{\it Ideals, Varieties, and Algorithms:
An Introduction to Computational Algebraic Geometry and Commutative Algebra,}
Springer Verlag, Heidelberg, 2007.

\bibitem{EKS}
G. Elekes, H. Kaplan and M. Sharir,
On lines, joints, and incidences in three dimensions,
{\it J. Combinat. Theory, Ser. A} 118 (2011), 962--977.
Also in arXiv:0905.1583.
%
\bibitem{Er46}
P.\ Erd{\H o}s,
On sets of distances of $n$ points,
{\it Amer. Math. Monthly} 53 (1946), 248--250.

\bibitem{Gu14}
L. Guth,
Distinct distance estimates and low-degree polynomial partitioning,
in arXiv:1404.2321.

\bibitem{GK}
L. Guth and N. H. Katz,
Algebraic methods in discrete analogs of the Kakeya problem,
{\it Advances Math.} 225 (2010), 2828--2839.
Also in arXiv:0812.1043v1.

\bibitem{GK2}
L. Guth and N. H. Katz,
On the Erd\H{o}s distinct distances problem in the plane,
{\it Annals Math.} 181 (2015), 155--190.
Also in arXiv:1011.4105.

\bibitem{Har}
J.\ Harris,
{\it Algebraic Geometry: A First Course},
Vol. 133. Springer-Verlag, New York, 1992.

\bibitem{Hart83}
R.\ Harshorne,
{\it Algebraic Geometry},
Springer-Verlag, New York. 1983.

\bibitem{KMSS}
H. Kaplan, J. Matou\v{s}ek, Z. Safernov\'a and M. Sharir,
Unit distances in three dimensions,
{\it Combinat. Probab. Comput.} 21 (2012), 597--610.
Also in arXiv:1107.1077.

\bibitem{KMS}
H. Kaplan, J. Matou\v{s}ek and M. Sharir,
Simple proofs of classical theorems in discrete geometry via
the Guth--Katz polynomial partitioning technique,
{\it Discrete Comput. Geom.} 48 (2012), 499--517.
Also in arXiv:1102.5391.

\bibitem{Kol}
J. Koll\'ar,
Szemer\'edi--Trotter-type theorems in dimension 3,
in arXiv:1405.2243.

\bibitem{Land}
J.\ M.\ Landsberg,
is a linear space contained in a submanifold? On the number of derivatives needed to tell,
{\it J. Reine Angew. Math.} 508 (1999), 53--60.

\bibitem{PS}
J. Pach and M. Sharir,
Geometric incidences,
in {\it Towards a Theory of Geometric Graphs} (J. Pach, ed.),
{\it Contemporary Mathematics}, Vol. 342,
Amer. Math. Soc., Providence, RI, 2004, pp. 185--223.

\bibitem{RSZ}
O. Raz, M. Sharir, and F. De Zeeuw,
Polynomials vanishing on Cartesian products: The Elekes--Szab\'o Theorem revisited,
manuscript, 2014.

\bibitem{Schmid}
J. Schmid,
On the affine B\'ezout inequality,
{\it Manuscripta Mathematica} 88(1) (1995), 225--232.

\bibitem{SSS}
M. Sharir, A. Sheffer, and N. Solomon,
Incidences with curves in $\reals^d$,
manuscript, 2014.

\bibitem{SSZ}
M. Sharir, A. Sheffer, and J. Zahl,
Improved bounds for incidences between points and circles,
{\it Combinat. Probab. Comput.}, in press.
Also in {\it Proc. 29th ACM Symp. on Computational Geometry} (2013), 97--106,
and in arXiv:1208.0053.

\bibitem{surf-socg}
M. Sharir and N. Solomon,
Incidences between points and lines in $\reals^4$,
{\it Proc. 30th Annu. ACM Sympos. Comput. Geom.}, 2014, 189--197.

\bibitem{SS4d}
M. Sharir and N. Solomon,
Incidences between points and lines in four dimensions,
in arXiv:1411.0777.

\bibitem{SS3dvar}
M. Sharir and N. Solomon,
Incidences between points and lines on a two-dimensional variety,
manuscript, 2014.

\bibitem{SoTa}
J.\ Solymosi and T.\ Tao,
An incidence theorem in higher dimensions,
{\it Discrete Comput. Geom.} 48 (2012), 255--280.

\bibitem{Sz}
L.~Sz{\'e}kely,
Crossing numbers and hard Erd{\H o}s problems in discrete geometry,
{\it Combinat. Probab. Comput.} {\bf 6} (1997), 353--358.

\bibitem{SzT}
E. Szemer\'edi and W.~T. Trotter,
Extremal problems in discrete geometry,
{\it Combinatorica} 3 (1983), 381--392.

\bibitem{T}
T. Tao,
From rotating needles to stability of waves: Emerging connections
between combinatorics, analysis, and PDE,
{\it Notices AMS} 48(3) (2001), 294--303.

\bibitem{tt:blog}
T. Tao,
The Cayley--Salmon theorem via classical differential geometry,
{\tt http://terrytao.wordpress.com}, March 2014.

\bibitem{w-lbanm-68}
H.~E. Warren,
Lower bound for approximation by nonlinear manifolds,
{\it Trans. Amer. Math. Soc.} 133 (1968), 167--178.

\bibitem{Za1}
J.\ Zahl,
An improved bound on the number of point-surface incidences in three dimensions,
\emph{Contrib. Discrete Math.} 8(1) (2013).
Also in arXiv:1104.4987.

\bibitem{Za2}
J.\ Zahl,
A Szemer\'edi-Trotter type theorem in $\reals^4$,
in arXiv:1203.4600.

\end{thebibliography}
\end{document}